\documentclass[11pt]{amsart}
\usepackage{amssymb}

\newcommand{\C}{\ensuremath{\mathbb{C}}}

\makeatletter
\newcommand{\sumprime}{\if@display\sideset{}{'}\sum%
            \else\sum'\fi}
\makeatother

\begin{document}

\numberwithin{equation}{section}

\newtheorem{theorem}{Theorem}[section]
\newtheorem{proposition}[theorem]{Proposition}
\newtheorem{conjecture}[theorem]{Conjecture}
\def\theconjecture{\unskip}
\newtheorem{corollary}[theorem]{Corollary}
\newtheorem{lemma}[theorem]{Lemma}
\newtheorem{observation}[theorem]{Observation}
\theoremstyle{definition}
\newtheorem{definition}{Definition}
\numberwithin{definition}{section}
\newtheorem{remark}{Remark}
\def\theremark{\unskip}
\newtheorem{question}{Question}
\def\thequestion{\unskip}
\newtheorem{example}{Example}
\def\theexample{\unskip}
\newtheorem{problem}{Problem}

\def\vvv{\ensuremath{\mid\!\mid\!\mid}}
\def\intprod{\mathbin{\lr54}}
\def\reals{{\mathbb R}}
\def\integers{{\mathbb Z}}
\def\N{{\mathbb N}}
\def\complex{{\mathbb C}\/}
\def\dist{\operatorname{dist}\,}
\def\spec{\operatorname{spec}\,}
\def\interior{\operatorname{int}\,}
\def\trace{\operatorname{tr}\,}
\def\cl{\operatorname{cl}\,}
\def\essspec{\operatorname{esspec}\,}
\def\range{\operatorname{\mathcal R}\,}
\def\kernel{\operatorname{\mathcal N}\,}
\def\linearspan{\operatorname{span}\,}
\def\lip{\operatorname{Lip}\,}
\def\sgn{\operatorname{sgn}\,}
\def\Z{ {\mathbb Z} }
\def\e{\varepsilon}
\def\p{\partial}
\def\rp{{ ^{-1} }}
\def\Re{\operatorname{Re\,} }
\def\Im{\operatorname{Im\,} }
\def\dbarb{\bar\partial_b}
\def\eps{\varepsilon}
\def\O{\Omega}
\def\Lip{\operatorname{Lip\,}}

\def\Hs{{\mathcal H}}
\def\E{{\mathcal E}}
\def\scriptu{{\mathcal U}}
\def\scriptr{{\mathcal R}}
\def\scripta{{\mathcal A}}
\def\scriptc{{\mathcal C}}
\def\scriptd{{\mathcal D}}
\def\scripti{{\mathcal I}}
\def\scriptk{{\mathcal K}}
\def\scripth{{\mathcal H}}
\def\scriptm{{\mathcal M}}
\def\scriptn{{\mathcal N}}
\def\scripte{{\mathcal E}}
\def\scriptt{{\mathcal T}}
\def\scriptr{{\mathcal R}}
\def\scripts{{\mathcal S}}
\def\scriptb{{\mathcal B}}
\def\scriptf{{\mathcal F}}
\def\scriptg{{\mathcal G}}
\def\scriptl{{\mathcal L}}
\def\scripto{{\mathfrak o}}
\def\scriptv{{\mathcal V}}
\def\frakg{{\mathfrak g}}
\def\frakG{{\mathfrak G}}

\def\ov{\overline}

\author{Siqi Fu}
\thanks
{The author was supported in part by NSF grant DMS
0805852.}
\address{Department of Mathematical Sciences,
Rutgers University-Camden, Camden, NJ 08102}
\email{sfu@camden.rutgers.edu}
\title[]  
{The Kobayashi metric in the normal direction \\
and the mapping problem}
\maketitle

\begin{abstract}
Estimates of the Kobayashi metric in the normal direction are used to study
the mapping problem in several complex variables.
\end{abstract}



\

\section{Introduction}\label{sec:intro}
Does every biholomorphic map between two smooth bounded
domains in $\C^n$ extend smoothly to the boundary?
This problem, known as the mapping problem, has played an
important role in several complex variables. (See the survey \cite{Forstneric93}
and the monograph \cite{BER99} for extensive coverage of this and related problems.)
Since the early 1970's,
biholomorphically invariant metrics, such as the Carath\'{e}odory, Bergman, and Kobayashi metrics, have been employed to study the mapping problem.
The Carath\'{e}odory metric was used by Vormoor \cite{Vormoor73} to show
that biholomorphic maps between strictly pseudoconvex domains with $C^2$-smooth boundaries extend to homeomorphisms on the closures of the domains, and by Henkin and Pinchuk to show that these biholomorphic maps \cite{Henkin73}, in fact, even proper holomorphic maps \cite{Pinchuk74}, are H\"{o}lder continuous of order $1/2$.  In a seminal paper~\cite{Fefferman74}, using the Bergman metric, Fefferman showed that biholomorphic maps between smooth bounded strictly pseudoconvex domains extend to diffeomorphisms on the closures.

The Kobayashi metric is the largest among the invariant metrics that are identical to the Poincar\'{e} metric on the unit disc and have the contracting property. It naturally becomes a useful tool to study the mapping problem.  Diederich and Forn{\ae}ss \cite{DiederichFornaess79} used it to establish the H\"{o}lder continuity of proper holomorphic maps onto bounded pseudoconvex domains with real analytic boundaries. The key step in this approach is to show that the Kobayashi metric blows up in all directions at a rate no less than $1/\delta^\eps$, where $\delta$ is the distance to the boundary from the point where the metric is evaluated and $\eps$ a positive constant. Once this estimate is established, the H\"{o}lder continuity of the proper holomorphic maps then follows from a classical method of Hardy and Littlewood,  and the Diederich-Forn{\ae}ss construction of bounded plurisubharmonic exhaustion functions \cite{DiederichFornaess77}. Boundary behavior of the Carath\'{e}odory and Kobayashi metrics on smooth bounded strictly pseudoconvex domains in $\C^n$ was studied by Graham in \cite{Graham75}. (See \cite{Fu95} for the Fefferman type asymptotic expansions of the metrics on these domains.) Estimates for the Kobayashi metric, as well as the Carath\'{e}odory and Bergman metrics, in terms of big constants and small constants, were obtained by Catlin \cite{Catlin89} for smooth bounded pseudoconvex domains of finite type in $\C^2$. (See \cite{Mcneal01} and references therein for related results.) Cho \cite{Cho92} extended the above-mentioned result of Diederich and Forn{\ae}ss to smooth bounded pseudoconvex domains of D'Angelo finite type in $\C^n$.  However, such estimates are not expected to hold for domains that are not of finite type. In fact, it is well known that for a smooth bounded domains in $\C^2$,  the $1/\delta^\eps$-estimate in all directions necessitates that the domain is pseudoconvex and of  finite type $\le 1/\eps$ (see, for example, \cite{Fu94b}).

Krantz \cite{Krantz92} showed that for any smooth bounded domain in $\C^n$, the Kobayashi metric blows up in the complex normal direction at a rate no less than $1/\delta^{3/4}$. (In sharp contrast, the Carath\'{e}odory and Bergman metrics can remain bounded in all directions due to the Hartogs phenomenon.) This result was later generalized in \cite{Fu94a, Fu94b} to estimates of the Kobayashi metric from below in
any direction, albeit the estimates depending only on the complex normal component of the direction. In particular, it was observed in \cite{Fu94a}  that a smooth bounded domain is pseudoconvex if and only if the Kobayashi metric satisfies the $1/\delta^\alpha$-estimate in the complex normal direction for some $\alpha>3/4$. These estimates were then used to study the mapping problem in \cite{Fu94b}.  In this expository paper, we elaborate and expand upon this work. This paper would not have been possible without the kind encouragement and generosity of Steve Krantz, to whom we are greatly indebted. We also thank Lina Lee for constructive suggestions.

\bigskip

\section{The Kobayashi metric in the normal direction}\label{sec:kobayashi}
\smallskip\smallskip

Let $\Omega$ be a bounded domain in $\C^n$.  Let $H(D, \Omega)$ be family
of holomorphic maps from the unit disc $D$ into $\Omega$. The Kobayashi
metric\footnote{Kobayashi defined the distance named after him by using chains of holomorphic discs (see \cite{Kobayashi70}). Royden \cite{Royden71} introduced  this infinitesimal metric and proved that the distance induced by it is identical to the Kobayashi distance.} on $\Omega$ is given by
\[
F_{\Omega}^{K}(z, X) = \inf \left\{ 1/\lambda  \mid
\exists f\in H(D, \Omega), f(0)=z, f'(0)=\lambda X, \lambda > 0 \right\},
\]
for $z\in \Omega$ and $X\in T^{1, 0}_z(\Omega)$. (In what follows, we
identify $X$ with a vector in $\C^n$.)
For $z, w\in\Omega$, let
\[
\ell_\Omega (z, w)=\inf\{\rho(a, b) \mid \exists f\in H(D, \Omega), f(a)=z, f(b)=w\},
\]
where $\rho$ is the Poincar\'{e} distance on $D$. We refer the reader to \cite{Krantz02} and \cite{JarnickiPflug93} for background material and extensive treatments on the Kobayashi and other invariant metrics.

The following localization property of the Kobayashi metric is due to Royden (\cite{Royden71}; see also \cite{Graham75}):

\begin{lemma}\label{lm:Royden} Let $U$ be any subdomain of $\Omega$. Then
\begin{equation}\label{eq:localization}
F^K_\Omega(z, X)\le F^K_U (z, X)\le \coth(\ell_\Omega(z, \Omega\setminus U))F^k_\Omega(z, X),
\end{equation}
for any $z\in U$ and $X\in \C^n$, where $\ell_\Omega(z, \Omega\setminus U)=\inf_{w\in\Omega\setminus U} \ell_\Omega(z, w)$.
\end{lemma}

The following lemma is well-known.

\begin{lemma}\label{lm:kob1}.
Let $\eps_0$ be a positive constant.  Suppose
that $\Omega$ is a domain in $\C$ such that, for each $z_0\in b\Omega$,  the
connected component of $\C\setminus \Omega$ containing $z_0$ also
contains a segment $\overline{z_0z_1}$ with length $\ge\eps_0$.  Then
\begin{equation}\label{eq:kob1}
{|X|\over d(z)}\ge F^K_{\Omega}(z, X) \ge \frac{1}{8}{|X|\over d(z)}
\end{equation}
for all $z\in \O$ with $d(z)\le \eps_0$, where $d(z)$ is the Euclidean
distance from $z$ to the boundary $b\Omega$.
\end{lemma}

\begin{proof} The first inequality holds for all $z\in\Omega$ and is easily
seen by comparing with the Kobayashi metric on the disc with center $z$
and radius $d(z)$.  We now prove the second inequality.

By homogeneity of the Kobayashi metric, we may assume that $X=1$.  Let $z$ be a point in $\O$ with $d(z)\le\eps_0$.  Let $z_0$ be a point in $b\Omega$ such that
$d(z)=|z_0-z|$.  Let $z_1$ be a point in the connected component of $\C\setminus\O$ containing $z_0$ such that $|z_1-z_0|=\eps_0$.
Let $\O_1=\C\setminus \{\overline{z_0z_1}\}$.  Then $\O_1\supset \O$.
After a translation and a rotation, we may assume without loss of generality that $z_0=0$ and $z_1=-\eps_0$. Let
$$
f(\zeta)=\left({\zeta\over \zeta+\eps_0}\right)^{1/2},
$$
where the square root is the branch obtained by deleting the negative
real axis. Then $f(\zeta)$ maps $\Omega_1$ into the right-half plane $\C^+$. It then follows from the length decreasing property of the Kobayashi metric that
\[
F^{K}_{\O}(z, 1)\ge F^{K}_{\C^+}(f(z), f'(z))= \frac{|f'(z)|}{2|\Re f(z)|}
          \ge \frac{|f'(z)|}{2|f(z)|}.
\]
A simple calculation then yields the desired inequality.
\end{proof}

We will write $z=(z_1, \ldots, z_n)=(\hat z, z_n)$ and use $B(z, R)$
to denote the ball in $\C^n$ with center $z$ and radius $R$.
For $0<k<1$, let $\Lambda(k)=\{z\in\C^n \mid -\Re z_n > k|z|\}$ be the cone
with axis in the negative $\Re z_n$-direction. Throughout the paper, we will use $C$ to denote constants which may be different in different appearances. The following proposition generalizes a result of Krantz \cite{Krantz92}. The proof is similar to
that in \cite{Krantz92}. The difference here is that whereas Krantz reduces the problem to an annulus in $\C$, we reduce it to the complement of a line segment in $\C$ (see Lemma~\ref{lm:kob1} above). In this way, we are able to estimate from below the Kobayashi metric in any direction, even though only the complex normal component of the direction contributes to the estimate.

\begin{proposition}\label{prop:ball1} Given constants $m\ge 2$, $R>0$, $A>0$, $k>0$, and $K>0$. Let
\[
\Omega_m=\{z\in B(0, R) \mid \Re z_n <A(|\hat z|^m+|z_n||z|)\}.
\]
Then there exist positive constants $C_1$ and $R_1$ such that
\begin{equation}\label{eq:kob2}
F^K_{\Omega_m}(z, X)\ge
C_1 \frac{|X_n|}{d^{1-\frac{1}{m}}_{\Omega_m}(z)},
\end{equation}
for all $z\in B(0, R_1)\cap \Lambda(k)$ and all $X\in\C^n$. Furthermore,
there exist positive constants $C_2$ and $C_3$ such that
\begin{equation}\label{eq:kob3}
C_2 \frac{|X_n|}{d^{1-\frac{1}{2m}}_{\Omega_m}(z)}\ge F^K_{\Omega_m}(z, X)\ge
C_3 \frac{|X_n|}{d^{1-\frac{1}{2m}}_{\Omega_m}(z)}.
\end{equation}
for all $z\in B(0, R_1)\cap \Lambda(k)$ and all $X\in\C^n$ with $|X|\le K|X_n|$.
\end{proposition}

\begin{proof}  Inequality \eqref{eq:kob2} and the second inequality
in \eqref{eq:kob3} have been proved in \cite{Fu94a}.  We include the
proofs here for completeness.

Assume that $R_1<R$. Let $z\in \Lambda(k)\cap B(0, R_1)\subset \Omega_m$ and $X\in\C^n$.   Let $\Phi(\zeta)=(\widehat{\Phi} (\zeta), \Phi_{n}
(\zeta))\colon D\to\Omega_m$ be a holomorphic map such that
\[
\Phi (0)=z,\quad
\Phi'(0)=\lambda X, \ \lambda>0.
\]
It follows from the  Cauchy integral formula that for $|\zeta|<1/2$,
\begin{equation}\label{eq:p1}
|\Phi(\zeta) -z|\le 2R |\zeta |
\end{equation}
and
\begin{equation}\label{eq:p2}
 |\Phi(\zeta)-z -\lambda
\zeta X |\le 4R |\zeta|^2.
\end{equation}
Denote $\delta =-\Re z_n$. From \eqref{eq:p1} and the assumption that $k|z|<\delta$, we obtain
\[
|\Phi(\zeta)|\le |z| + 2R|\zeta |\le (1/ k) \delta + 2Rc\delta^{1/
m}\le (1/k+2Rc) \delta^{1/m}
\]
for  $|\zeta|< c\delta^{1/m}$, provided $\delta<1$. (Note that $\delta<R_1$ and we can take $R_1$ to be sufficiently small.) Therefore,
\begin{align}
\Re \Phi_n (\zeta) &< A\left( |\widehat{\Phi} (\zeta)|^m + |\Phi_n
(\zeta)|\cdot |\Phi (\zeta)|\right)\label{eq:p4}\\
&\le \frac{1}{2}\left(\delta+ |\Phi_n(\zeta)|\right),\notag
\end{align}
when $\delta$ and $c$ are sufficiently small. Let $f(\zeta)=\Phi_n(c\delta^{1/m}\zeta)$. Then $f(0)=-\delta$,
$f'(0)=c\delta^{1/m}\lambda X_n$, and
\[
f(D)\subset \C\setminus\{w\in\C \mid \Im w=0, \Re w\ge\delta\}.
\]
It then follows from Lemma~\ref{lm:kob1} that
\begin{equation}\label{eq:p5}
c\delta^{1/m}\lambda |X_n| \le C\delta,
\end{equation}
from which \eqref{eq:kob2} follows.

We now prove the second inequality in \eqref{eq:kob3}.
From \eqref{eq:p2}, we have
\begin{align}
 |\Phi (\zeta)|&\le  |z| + (\lambda |X|)|\zeta | + 4R|\zeta|^2\notag\\
&\le (1/k)\delta + (\lambda |X|)
c\delta^{1/2m} +
4c^{2} R\delta^{1/m}\notag\\
&\le c^{1/2}\left(\delta^{1/m} + (\lambda
|X|)\delta^{1/2m}\right),\label{eq:p6}
\end{align}
for $|\zeta |<c\delta^{1/2m}$, provided $c$ and $\delta$ are sufficiently small. From \eqref{eq:p4} and \eqref{eq:p6}, we then have
$$
\Re \Phi_{n}(\zeta)< \frac{\delta + (\lambda |X|)^m \delta^{1/2}}{2} +
\frac{1}{2} |\Phi_n (\zeta)|
$$
for $|\zeta|<c\delta^{1/2m}$, when $c$, $\delta$  are sufficiently small.
Let $g(\zeta)=\Phi_n(c\delta^{1/2m}\zeta)$. Then $g(0)=-\delta$,
$g'(0)=c\delta^{1/2m}\lambda X_n$, and
$$
g(D)\subset\C\setminus \left\{ w\in
\C\enspace\bigm|\enspace
\Im w=0, \Re w\ge \delta +(\lambda |X|)^{m} \delta ^{1/2} \right\}.
$$
It follows from Lemma~\ref{lm:kob1} that
\[
c\delta^{1/2m}\lambda |X_n|\le C\left(\delta+ (\lambda |X|)^{m} \delta ^{1/2}\right).
\]
From \eqref{eq:p5} and the assumption $|X|\le K|X_n|$, we then
have
\[
c\delta^{1/2m}\lambda |X_n|\le C\left(\delta+ (K\lambda |X_n|)^{m} \delta ^{1/2}\right)\le C(\delta+\delta^{m-\frac{1}{2}})\le C\delta.
\]
We thus conclude the proof of the second inequality in \eqref{eq:kob3}.

We will actually prove a slightly stronger version of the first inequality
in \eqref{eq:kob3} in the next proposition.
\end{proof}

\begin{proposition}\label{prop:ball2} Given constants $m\ge 2$, $R>0$, $A>0$, $B>0$, $k>0$, and $K>0$. Let
\[
\widetilde\Omega_m=\{z\in B(0, R) \mid r(z)=\Re z_n -A|z_1|^m+B(\sum_{j=2}^n |z_j|^m +|z_n||z|)<0\}.
\]
Then there exist positive constants $R_1$ and $C$ such that
\begin{equation}\label{eq:kob4}
F^K_{\widetilde\Omega_m}(z, X)\le C \frac{|X_n|}{d^{1-\frac{1}{2m}}_{\widetilde\Omega_m}(z)}
\end{equation}
for all $z\in B(0, R_1)\cap \Lambda(k)$ and all $X\in\C^n$ with
$|X|\le K|X_n|$.
\end{proposition}

\begin{proof} By homogeneity of the Kobayashi metric,  we may assume that $X_n=1$.  Let $z=(z_1,z_2,\ldots,
z_n)\in\Lambda(k)\cap\widetilde\O_m$ and let $\delta =-\Re z_n$.  Define
$$ \Phi_{\delta}(\zeta) =(\Phi_{1\delta}(\zeta),
\Phi_{2\delta}(\zeta), \ldots, \Phi_{n\delta}(\zeta))
$$
by
\begin{align*}\Phi_{1\delta} &= z_1 +c\delta^{1-1/2m}X_{1}\zeta
+b^{-1}\zeta^2;\cr
\intertext{and}
\Phi_{k\delta} &=  z_{k}+c\delta^{1-1/2m}X_{k}\zeta,
\quad\hbox{\rm for}\quad  2\le k\le n,\cr
\end{align*}
where $b$ and $c$ are sufficiently small constants to be chosen. Then $\Phi_{\delta}(0)=z$,  $\Phi'_{\delta}(0)=c\delta^{1-1/2m}X$.
It suffices to prove that $\Phi_{\delta}(D)\subset \widetilde\Omega_m$ for
sufficiently small $c$ and $\delta$.

If $1>|\zeta|\ge\delta^{1/2m}$, we have
\begin{align*}
|\Phi_{1\delta}(\zeta)|&\ge b^{-1}|\zeta|^2-(k^{-1}\delta+c\delta^{1-1/2m} |X_1|\|\zeta|)\\
&\ge b^{-1}|\zeta|^2-(k^{-1}+cK)|\zeta|^{2m}\ge |\zeta|^2
\end{align*}
when $b$ is sufficiently small and $c<1$.  Thus
\begin{align*}
r(\Phi_\delta(\zeta))&\le-\delta+c\delta^{1-1/2m}\Re\zeta -A|\zeta|^{2m}-C\delta^{3/2}\\
&\le -\delta/2-(A-c)|\zeta|^{2m}<0,
\end{align*}
when $c$ and $\delta$ are sufficiently small.

If $|\zeta|<\delta^{1/2m}$, then it is easy to see that
\[
r(\Phi_\delta(\zeta))\le -\delta+c\delta^{1-1/2m}\Re\zeta+c^{1/2}\delta<0,
\]
for sufficiently small $c$ and $\delta$.  We thus conclude the proof of
the proposition.
\end{proof}

We now use the above propositions to estimate the Kobayashi metric in the normal direction for a general domain in $\C^n$. Let $\Omega=\{z\in\C^n \mid r(z)<0\}$ be a bounded domain in $\C^n$ with $C^k$-smooth boundary. The defining function $r(z)$ is always chosen to be in the same smoothness class as the boundary and $dr\not=0$ on $b\Omega$.  One useful choice of the defining function is the signed
Euclidean distance function:
\[
\delta_\Omega(z)=
\begin{cases} -d(z, b\Omega),   & \text{if $z\in\Omega$},\\
d(z, b\Omega), &\text{if $z\in\C^n\setminus\Omega$}.
\end{cases}
\]
We refer the reader to \cite{KrantzParks81} and references therein for
more information on the distance function. In particular, it was
shown there that when $b\Omega$ is of class $C^k$, $k\ge 2$, then $\delta_\Omega(z)$ is $C^k$-smooth in a neighborhood of $b\Omega$.
This, however, is not true when $k=1$. Nonetheless, if $b\Omega$ is of class $C^{1, 1}$, so is $\delta_\Omega(z)$ in a neighborhood of $b\Omega$.

\begin{theorem}\label{th:kob1}
Let $\Omega=\{z\in\C^n \mid r(z)<0\}$ be a bounded domain with $C^{1, 1}$-smooth boundary. Then there exists a constant $C>0$ such that
\begin{equation}\label{eq:n1}
F^K_\Omega(z, X)\ge C\frac{|\langle \partial r(z), X\rangle|}{
|r(z)|^{1/2}}
\end{equation}
for all $z\in\Omega$ and all $X\in \C^n$.
\end{theorem}

\begin{proof} This theorem has been proved in \cite{Fu94b} under the assumption
that $b\Omega$ is $C^\infty$-smooth. The same proof works when $b\Omega$ is $C^{1, 1}$-smooth, in light of the above discussion. We provide the detail below.

It is easy to see that \eqref{eq:n1} is independent of the choice of the defining function. It suffices to establish
\eqref{eq:n1} for $\delta_\Omega(z)$ in a neighborhood of $b\Omega$.
Let $U$ be a neighborhood of $b\Omega$ so that $\delta_\Omega\in C^{1, 1}$ and there is a projection $\pi\colon U\to b\Omega$ such that $|z-\pi(z)|=d(z, b\Omega)$.  Let $p\in U\cap\Omega$. After a translation and a unitary transformation, we may assume that $\pi(p)$ is the origin and the outward normal direction at $\pi(p)$ is the positive $\Re z_n$-axis. Hence there exist positive constants $R$ and $A$, independent of $p$, such that
\[
\Omega\cap B(0, R)\subset \{z\in B(0, R) \mid \Re z_n<A(|\hat z|^2+|z_n||z|)\}.
\]
It follows from Proposition~\ref{prop:ball1} that
\[
F^K_{\Omega\cap B(0, R)}(p, X)\ge C\frac{|X_n|}{|\delta_\Omega(p)|^{1/2}}.
\]
when $\delta_\Omega(p)$ is sufficiently small. (The constants can be chosen to be independent of $\pi(p)$.) Since in the new coordinates, $\partial\delta_\Omega(p)/\partial z_j =0$, $1\le j\le n-1$, and $\partial\delta_\Omega(p)/\partial z_n=1/2$. We then obtain \eqref{eq:n1} after applying Lemma~\ref{lm:Royden}.
\end{proof}

\begin{remark} The above theorem is sharp. For example, let $\Omega$ be a bounded domain in $\C^2$ locally defined near the origin by $r(z)=\Re z_2-|z_1|^2<0$. Let
$p_\delta=(-\delta, 0)$ and $X=(\delta^{-1/2}, 1)$, then
\[
F^K_\Omega(p_\delta, X)\approx \frac{1}{\delta^{1/2}},
\]
when $\delta$ is sufficiently small.
\end{remark}

\begin{theorem}\label{th:n2}
$(1)$ Let $\Omega=\{z\in\C^n \mid r(z)<0\}$ be a bounded domain with $C^2$-smooth boundary. If there exist constants $C>0$ and $\alpha>1/2$ such
that
\begin{equation}\label{eq:n2}
F^K_\Omega(z, X)\ge C\frac{|\langle \partial r(z), X\rangle|}{
|r(z)|^{\alpha}}
\end{equation}
for all $z\in\Omega$ and all $X\in \C^n$, then $\Omega$ is pseudoconvex.

$(2)$ Conversely, if $\Omega=\{z\in\C^n \mid r(z)<0\}$ is a bounded pseudoconvex domain with $C^3$-smooth boundary, then there exists a constant $C>0$ such that
   \begin{equation}\label{eq:n3}
F^K_\Omega(z, X)\ge C\frac{|\langle \partial r(z), X\rangle|}{
|r(z)|^{2/3}}
\end{equation}
for all $z\in\Omega$ and all $X\in \C^n$.
\end{theorem}

\begin{proof} Let $p^0\in b\Omega$. After a translation and a rotation, we can assume that $p^0$ is the origin and the outward normal direction at $p^0$ is the positive real $\Re z_n$-axis. Furthermore, after the following simple change of coordinates
\[
(\hat z, z_n)\mapsto (\hat z, z_n+\sum_{j, k=1}^n r_{z_j z_k}(0)z_j z_k),
\]
we may assume that
\[
r(z)=\Re z_n + \sum_{j, k=1}^n a_{jk}z_j\bar z_k +o(|z|^2).
\]
Notice all these changes of coordinates preserve the outward normal direction at $p^0$.

We now prove (1) by contradiction. Suppose that the Levi-form of $r(z)$
is not semi-positive at some $p^0\in b\Omega$.  After changes of coordinates as above, we assume that $p^0$ is the origin and
$$
r(z)=\Re z_n + \sum_{j,k=1}^{n} a_{jk}z_{j}\bar{z_{k}} + o(|z|^2)
$$
for $z$ near $p^0$. Since the Levi-form is not semi-positive definite
at $p^0$, the matrix
$\left({a_{jk}}\right)_{ 1\le i,j  \le
n- 1}$   has at least one
negative eigenvalue .  After a unitary transformation in
the $\hat z=(z_1, z_2, \ldots, z_{n-1})$ variables, we may assume
$$
r(z)=\Re z_n +\sum_{j=1}^{n-1} c_j |z_j|^2 + 2\Re\sum_{j=1}^{n-1} a_{jn}z_{j}\bar{z_{n}}+a_{nn}|z_n|^2 + o(|z|^2)
$$
with $c_1<0$. Write $c_1=-2a$.  Applying the inequality $|ab|\le \eps |a|^2+(1/\eps)b^2$, we find positive constants $R$ and $B$ such that
$$
r(z)\le \hat r(z)=\Re z_n -a|z_1|^2 + B \sum_{j=2}^n |z_j|^2<0,
$$
for $z\in B(0, R)$. Let $\widehat\Omega=\{z\in B(0, R) \mid \hat r(z)<0\}$.
Then $\Omega\cap B(0, R)\supset \widehat\Omega$. Let
\[
p_\delta =(0, \ldots, 0, -\delta), \quad X_\delta=(\frac{\delta^{-1/2}}{\sqrt{a}}, 0, \ldots, 0, \frac{1}{2}),\quad \Phi_\delta (\zeta)=(\frac{\zeta}{\sqrt{a}}, 0, \ldots, 0, -\delta+\frac{\delta^{1/2}\zeta}{2}).
\]
Then $\Phi_\delta(0)=p_\delta$ and $\Phi'_\delta(0)=\delta^{1/2}X_\delta$.
By splitting into two cases, $|\zeta|<\delta^{1/2}$ and $|\zeta|\ge\delta^{1/2}$, as in the proof of Proposition~\ref{prop:ball2}, we obtain that $\Phi_\delta(D)\subset \widehat\Omega$ for sufficiently small $\delta$.  Hence
\begin{equation}\label{eq:contra}
F^K_{\widehat\Omega}(P_\delta, X_\delta)\le C\delta^{-1/2}.
\end{equation}
Since
\[
\frac{\partial r}{\partial z_n}(p_\delta)=-\frac{1}{2}+O(\delta) \quad \text{ and } \quad
\frac{\partial r}{\partial z_j}(p_\delta)=O(\delta), \ 1\le j\le n-1,
\]
we have $|\langle \partial r(p_\delta), X_\delta\rangle |\le 1$ for sufficiently small $\delta$. Therefore, by \eqref{eq:n2} and  Lemma~\ref{lm:Royden}, we have
\[
\delta^{-\alpha}\le C F^K_\Omega(p_\delta, X_\delta)\le C F^K_{\Omega\cap B(0, R)}(p_\delta, X_\delta)\le C F^K_{\widehat\Omega}(P_\delta, X_\delta),
\]
which contradicts
\eqref{eq:contra}. We have thus proved the first part of the theorem.

The proof of the second part is similar to that of Proposition~\ref{prop:ball2}. Let $p\in\Omega$ be sufficiently close
to the boundary $b\Omega$. Let $p^0$ be the projection of $p$ onto
the  boundary such that $\delta=|p-p^0|$ is the distance from $p$ to the boundary . Proceeding as before, after changes of coordinates that preserves the outward normal direction
at $p^0$, we have $p^0=0$ and
$$
r(z)=\Re z_n + \sum_{j,k=1}^{n} a_{jk}z_{j}\bar{z_{k}} + O(|z|^3).
$$
Since $b\Omega$ is pseudoconvex, the hermitian matrix $\left({a_{jk}}\right)_{ 1\le i,j  \le
n- 1}$ is semi-positive definite.  Hence
\begin{align*}
r(z) &\ge \Re z_n+2\Re\sum_{j=1}^{n-1} a_{jn}z_j\bar z_n+a_{nn}|z_n|^2-C|z|^3\\
&\ge \Re z_n -A(|\hat z|^3+|z_n||z|),
\end{align*}
for some positive constants $C$ and $A$. It then follows from Proposition~\ref{prop:ball1} that
\[
F^K_\Omega(p, X)\ge C\frac{|X_n|}{|r(z)|^{2/3}}
\]
when $p$ is sufficiently close to $b\Omega$. Since
\[
\langle\partial r(p), X\rangle=X_n/2+O(\delta),
\]
and the constants in the above estimates can be chosen to be independent of $p^0$, we then obtain \eqref{eq:n3}. \end{proof}

\begin{remark} Inequality \eqref{eq:n3} is not supposed to be sharp: One should be able to replace $2/3$ by any positive constant less than 1. However, the examples of Krantz \cite{Krantz93} and Forn{\ae}ss and Lee \cite{FornaessLee08} show that one cannot
replace $2/3$ by 1. These estimates on the Kobayashi metric in the complex normal direction are closely related to completeness of the metric. For example, it is easy to see that if
\[
F^K_\Omega(z, X)\ge C \frac{\left|\langle\partial r(z), X\rangle\right|}{g(r(z))},
\]
where $g$ is a positive function with $\int_{-\infty}^0 (g(t))^{-1}\, dt=\infty$, then the Kobayashi metric is complete on $\Omega$.
\end{remark}

\section{the mapping problem}\label{sec:maps}

In this section, we explain how estimates of the Kobayashi metric in the normal
direction in the previous section can be used to study the mapping problem.
Let $\Omega$ be a bounded domain in $\C^n$ with $C^2$-smooth boundary. Let $\delta(z)=\delta_\Omega(z)$ be the signed Euclidean distance to the boundary $b\Omega$. For $p$ sufficiently close to $b\Omega$, denote by
\[
n_p=\sum_{j=1}^n \left(\frac{\partial \delta}{\partial x_j}\frac{\partial }{\partial x_j}+\frac{\partial \delta}{\partial y_j}\frac{\partial }{\partial y_j}\right)
\]
the outward (real) normal direction and by
\[
N_p=2\sqrt{2}\sum_{j=1}^n \frac{\partial\delta}{\partial \bar z_j}\frac{\partial}{\partial z_j}.
\]
the complex normal direction. Note that $n_p=\sqrt{2}\Re N_p$.

Let $\Phi$ be a proper holomorphic map from $\Omega_1$ onto $\Omega_2$; both domains are assumed to have $C^2$-smooth boundaries.  Let $n_j$ and $N_j$ be the real and complex normal directions respectively, defined as above, of the domains $\Omega_j$, $j=1, 2$. Write $\delta_j(z)=\delta_{\Omega_j}(z)$.

\begin{definition}
The map $\Phi$ is said to {\it uniformly preserve the real normal direction} if there exist a neighborhood $U$ of $b\Omega_1$ and a constant $C>0$ such that
\begin{equation}\label{eq:real}
\left\|\Phi_{*p}(n_{1p})\right\|\le C\left|\langle d\delta_2 (\Phi(p)), \Phi_{*p}(n_{1p})\rangle\right|
\end{equation}
for all $p\in \Omega\cap U$.  It is said to {\it uniformly preserve the complex normal
direction} if
there exist a neighborhood $U$ of $b\Omega_1$ and a constant $C>0$ such that
\begin{equation}\label{eq:complex}
\left\|\Phi_{*p}(N_{1p})\right\|\le C\left|\langle \partial\delta_2 (\Phi(p)), \Phi_{*p}(N_{1p})\rangle\right|
\end{equation}
for all $p\in \Omega\cap U$.
\end{definition}

Geometrically, \eqref{eq:real} says that
the image of the real normal direction of the level set of $\delta_1(z)$ at $p\in\Omega_1$ under $\Phi$ stays uniformly away from the real tangent space of the level set of $\delta_2(z)$ at $\Phi(p)$. Namely,
\[
|\measuredangle\langle n_{2\Phi(p)}, \Phi_{*p}(n_{1p})\rangle-\pi/2|\ge c>0,
\]
for some constant $c>0$, where $\measuredangle\langle \cdot, \ \cdot\rangle$ denotes the angle
between the two vectors.  Similarly, \eqref{eq:complex} says that the image of the complex normal direction of the level set of $\delta_1(z)$ at $p_1\in\Omega_1$ under $\Phi$ stays away from the complex tangent space of the level set of $\delta_2(z)$ at $\Phi(p_1)$. Namely,
\[
|\measuredangle\langle N_{2\Phi(p)}, \Phi_{*p}(N_{1p})\rangle-\pi/2|\ge c>0,
\]
for some constant $c>0$. Obviously, any proper holomorphic map between domains on the plane uniformly preserves the complex normal direction. Also, if $\Phi$ is a $C^1$-diffeomorphism on the closures, then it uniformly preserves the real normal direction. It is easy to see that \eqref{eq:real} implies \eqref{eq:complex}
and that \eqref{eq:real} is equivalent to
\begin{equation}\label{eq:normal}
\left\|\frac{\partial\Phi}{\partial n_1}(p)\right\|\le C \left| \frac{\partial \delta_2\circ\Phi}{\partial n_1}(p)\right|.
\end{equation}
Note that if $\delta_2$ is plurisubharmonic, then by the Hopf lemma,
$\partial \delta_2\circ \Phi /\partial n_1\ge C>0$. Also, the definitions
are independent of the choices of defining functions and changes of coordinates.

Suppose that $\Omega_1$ is pseudoconvex. Then it follows from the Diederich-Forn{\ae}ss \cite{DiederichFornaess77} construction of bounded plurisubharmonic exhaustion functions on $\Omega_1$ and the Hopf lemma that $\Phi$ satisfies the following property (see \cite{DiederichFornaess79}): There exists constants $C>0$ and $\alpha \in (0, \, 1]$ such that
\begin{equation}\label{eq:alpha}
d_2(\Phi(z))\le C d_1^\alpha(z) \tag{${\rm DF}_\alpha$}
\end{equation}
for all $z\in\Omega_1$, where $d_k$ is the Euclidean distance to $b\Omega_k$, $k=1, 2$.

Denote by $\Lip_\alpha(\Omega)$ the standard Lipschitz class of order $\alpha$. The following two theorems were proved in \cite{Fu94b}.

\begin{theorem}\label{th:smooth}
Let $\Phi\colon\Omega_1\to\Omega_2$ be a proper holomorphic map between two bounded
domains in $\C^n$ with $C^2$-smooth boundary. If $\Phi$ uniformly preserves
the complex normal direction and satisfies property \eqref{eq:alpha} for
some $\alpha\in (0, \, 1]$. Then $\Phi\in\Lip_{\frac{1}{2}\alpha}(\Omega_1)$.
\end{theorem}

\begin{proof} Let $U$ be a tubular neighborhood of $b\O_1$  such that the orthogonal projection $\pi\colon U\to b\O_1$  is a smooth retraction.  Choose $\eps_0>0$ sufficiently small such that
$$
b\O_\eps =\left\{ z\in \O; \enspace d(z)=\eps \,\right\}\subset U
$$
for all $\eps\in (0, \eps_0)$   and $\pi \colon b\O_{\eps}\to b\O$ is a
diffeomorphism.  Let $p\in b\O_1$ and let $n_{1p}$ and $N_{1p}$ be the real and complex normal
directions of $b\O_1$ at $p$ respectively.  Set $p(t)=p-tn_{1p}$.
For $t_1, t_2\in (0, \eps_0)$ with $t_1\le t_2$,
\begin{equation}\label{eq:c1}
\left\| \Phi(p(t_1))-\Phi (p(t_2))\right\| =\left\| \int_{t_1}^{t_2} {d\Phi
(p(t))\over dt}\, dt \right\|
\le \int_{t_1}^{t_2} \left\| \Phi_{*p(t)} (n_{p(t)})\right\| \, dt.
\end{equation}
Since $\Phi$ uniformly preserves the complex normal direction,
\begin{equation}\label{eq:c2}
\left\| \Phi_{*p(t)} (n_{1p(t)})\right\|=\left\|\Phi_{*p(t)}(N_{1p(t)})\right\|
                 \le C\left|\langle \Phi_{*p(t)}(N_{1p(t)}), \enspace \partial \delta_2(\Phi(p(t)))\rangle\right|.
\end{equation}
By Theorem~\ref{th:kob1} and the length decreasing property of the Kobayashi metric, we have
\begin{equation}\label{eq:decreasing}
{|X|\over d_1(z)}\ge F^K_{\O_1}(z, X)\ge F^K_{\O_2} (\Phi (z), \Phi_{*z}(X))\ge C { \left|\langle \Phi_{*z}(X),   \partial
d_2(\Phi(z))\rangle\right|\over
d_2^{1/2}(\Phi(z))}.
\end{equation}
It then follows from property \eqref{eq:alpha} and \eqref{eq:decreasing} that
\begin{equation}\label{eq:c3}
\left|\langle \Phi_{*z}(X),   \partial
d_2(\Phi(z))\rangle\right|\le C {|X|\over d_1^{1-{1\over
2}\alpha}(z)}.
\end{equation}
Combining \eqref{eq:c1}, \eqref{eq:c2}, and \eqref{eq:c3}, we obtain
$$
\left\| \Phi(p(t_1))-\Phi(p(t_2))\right\|\le C\int_{t_1}^{t_2}{1\over
t^{1-{1\over 2}\alpha}}\, dt\le C (t^{{1\over 2}\alpha}_2 -t^{{1\over 2}\alpha}_1)\le C (t_2 -t_1)^{{1\over 2}\alpha}.
$$
It then follows from a result of Krantz (\cite{Krantz80}, Theorem 2.2) that $\Phi\in \Lip_{{1\over 2}\alpha}(\O_1)$.
\end{proof}

\begin{theorem}\label{th:map2}
Let $\Phi\colon\Omega_1\to\Omega_2$ be a proper holomorphic map between two bounded pseudoconvex domains with $C^3$-smooth boundary in $\C^n$. If $\Phi$ uniformly preserves the complex normal direction, then $\Phi\in\Lip_\beta(\Omega)$ for any $\beta<2/3$.
\end{theorem}

\begin{proof}  As we have discussed above, $\Phi$ satisfies property \eqref{eq:alpha}
for some $\alpha>0$.  Using Theorem~\ref{th:n2} (2) and following the proof of Theorem~\ref{th:smooth},  we have $\Phi\in\Lip_{\frac23\alpha}(\Omega)$. Once we know that $\Phi$ is continuous, we can localize the problem and choose $\alpha$ to be arbitrarily close to 1 (see pp. 591 in \cite{DiederichFornaess79}).
\end{proof}

Recall that a smooth bounded domain $\Omega$ is said to satisfy {\it condition R} if the Bergman projection $P$ maps $C^\infty(\overline{\Omega})$ into itself. A theorem of Bell and Ligocka \cite{BellLigocka80} says that biholomorphic maps between smooth bounded domains satisfying condition R extend smoothly to the boundaries. (See \cite{Bell90} for an exposition in this direction.) Barrett \cite{Barrett84} constructed a smooth bounded domain $\Omega$ such that $P(C^\infty_0(\Omega))$ is not even contained in $L^{2+\eps}(\Omega)$ for any $\eps>0$.  Nonetheless, it follows from a result of Barrett \cite{Barrett86} that for a strictly starlike smooth bounded domain $\Omega$ in $\C^n$,
\begin{equation}\label{eq:barrett}
P(C^\infty_0(\Omega))\subset L^{2+\eps}(\Omega) \tag{${\rm B}_\eps$},
\end{equation}
for some $\eps>0$. It is well-known that all smooth bounded pseudoconvex domains satisfy the above property (\cite{Kohn99}; see also \cite{BoasStraube99, BerndtssonCharpentier00}).  Barrett proved that for a given Diederich-Forn{\ae}ss \cite{DiederichFornaess77b} worm domain $\Omega$, the Bergman projection does not preserve the $L^2$-Sobolev spaces $W^k(\Omega)$ for sufficiently large $k$ \cite{Barrett92}.  M.~Christ resolved a long standing conjecture by showing that the Diederich-Forn{\ae}ss worm domains do not satisfy condition R (\cite{Christ96}). We refer the reader to the surveys \cite{Christ99, BoasStraube99, DangeloKohn99}, the book \cite{ChenShaw99}, and the recent paper of Straube \cite{Straube08} for more information on the related regularity theory in the $\bar\partial$-Neumann problem.

Let $\Phi\colon\Omega_1\to\Omega_2$ be a biholomorphic map between two smooth bounded domains in $\C^n$.  Suppose that the domain $\Omega_1$ satisfies property \eqref{eq:barrett} and  $\Omega_2$ satisfies the following property: There exists a function $g\in C^\infty_0(\Omega_2)$ such that
\begin{equation}\label{eq:bell}
|P_2(g)|\ge c>0
\end{equation}
for some $c>0$, where $P_2$ is the Bergman projection of $\Omega_2$. Lempert \cite{Lempert86} showed that under these conditions, $\Phi$ satisfies property \eqref{eq:alpha} with
\begin{equation}\label{eq:lempert}
\alpha=\frac{\eps}{2n(2+\eps)}.
\end{equation}
When $\Omega_2$ has real analytic boundary, then there exists a function $g\in C^\infty_0(\Omega_2)$ such that $P_2(g)=1$ (\cite{Bell81}). Hence property \eqref{eq:bell} is satisfied. Lempert further proved that in this case, $\Phi\in\Lip_\beta(\Omega_1)$ for some $\beta>0$ (\cite{Lempert86}). The following corollary exhibits a connection between the mapping problem and the problem of Lu Qi-Keng on the vanishing of the Bergman kernel. (See \cite{BFS99} for more information about Lu Qi-Keng's problem.)

\begin{corollary}\label{cor:lu}
Let $\Phi\colon\Omega_1\to\Omega_2$ be a biholomoprhic map between two smooth bounded domains in $\C^n$. Suppose that $\Omega_1$ satisfies property \eqref{eq:barrett} and
$\Omega_2$ satisfies the following property: There exists a point $w_0\in\Omega_2$ and
a constant $c>0$ such that $|K_2(w, w_0)|\ge c$ for all $w\in\Omega_2$. If $\Phi$ uniformly preserves the complex normal direction, then $\Phi\in\Lip_{\frac{1}{2}\alpha}(\Omega_1)$, where $\alpha>0$ is the constant given by \eqref{eq:lempert}.
\end{corollary}

\begin{proof} Let $\theta$ be a radially symmetric smooth function, compactly supported in the unit ball, such that its integral is 1.  Let $\theta_{w_0}(w)=t^{-2n}\theta((w-w_0)/t)$ where $t=d_2(w_0)$. Then $|P_2(\theta_{w_0})(w)|=|K_2(w, w_0)|\ge c>0$. Therefore, \eqref{eq:bell} and hence \eqref{eq:alpha} are satisfied, by the above-mentioned result of Lempert. The corollary then follows from Theorem~\ref{th:smooth}.
\end{proof}

\bibliography{survey}
\providecommand{\bysame}{\leavevmode\hbox
to3em{\hrulefill}\thinspace}

\end{document}